\documentstyle[12pt]{article}
\def\Bbb R{{\rm \bf R}}
\def\proclaim#1{\vskip2mm{\bf #1}\em}
\def\endproclaim{\em \vskip2mm}
\def\tag#1{\eqno(#1)}
\def\gathered{\begin{array}{c}}
\def\endgathered{\end{array}}
\def\text{\mbox}
\begin{document}

\title {On uniqueness in the inverse obstacle problem
via the positive supersolutions of the Helmholtz equation}
\author{Masaru IKEHATA\footnote{
Department of Mathematics,
Graduate School of Engineering,
Gunma University, Kiryu 376-8515, JAPAN}
}
\date{Final, February 10, 2012}
\maketitle
\begin{abstract}
This paper is concerned with an inverse obstacle scattering
problem of an acoustic wave for a single incident plane wave and a wave
number.  The Colton-Sleeman theorem states the unique recovery
of sound-soft obstacles with a smooth boundary from the far-field pattern
of the scattered wave for a single incident plane wave at a fixed wave number.
The wave number has a bound given by the first Dirichlet eigenvalue of the negative Laplacian
in an open ball that contains the obstacles.
In this paper, another proof of the Colton-Sleeman theorem
that works also for the case when we have a known {\it unbounded} set that contains obstacles
is given.
Unlike original one, the proof given here is not based on the monotonicity of
the first Dirichlet eigenvalue of the negative Laplacian.
Instead, it relies on a {\it positive supersolution} of the Helmholtz equation
in a known domain that contains obstacles.
Some corollaries which are new and not covered by the Colton-Sleeman
Theorem are also given.

\noindent
AMS: 35R30


\end{abstract}


\section{Introduction}

This paper is concerned with an inverse obstacle scattering
problem of acoustic wave for a single incident plane wave and wave
number.

Let $D\subset\Bbb R^m$, $m=2,3$ be a bounded domain with $C^2$ boundary such that
$\Bbb R^m\setminus\overline D$ is connected.
The total wave field $u$ outside {\it sound-soft} obstacle $D$ takes the form $u(x;d,k)=e^{ikx\cdot d}+w(x)$
with $k>0$, $d\in S^{m-1}$
and satisfies
$$\begin{array}{c}
\displaystyle
\triangle u+k^2 u=0\,\,\text{in}\,\Bbb R^m\setminus\overline D,\\
\\
\displaystyle
u=0\,\,\text{on}\,\partial D,\\
\\
\displaystyle
\lim_{r\longrightarrow\infty}r^{(m-1)/2}\left(\frac{\partial w}{\partial r}-ikw\right)=0.
\end{array}
$$
The last condition above is called the Sommerfeld radiation condition \cite{CK}.

The scattered wave $w=u-e^{ikx\cdot d}$ for fixed $d$ and $k$ has
the asymptotic behaviour
$$\displaystyle
w(r\varphi)=\frac{\displaystyle
e^{ikr}}{\displaystyle
r^{(m-1)/2}}F_D(\varphi,d;k)+O\left(\frac{1}{r^{(m+1)/2}}\right),\,\, r\longrightarrow\infty
$$
uniformly for $\varphi\in S^{m-1}$, and coefficient
$F_D(\varphi,d;k)$ is called the {\it far-field pattern}.

We consider the {\it uniqueness issue} of the inverse problem:
determine $D$ from $F_D(\,\cdot\,,d;k)$ for a fixed $d$ and $k$.

This is a well-known open problem, and the complete answer is yet
unknown \cite{IS}. However, there is a partial result with a constraint on
the range of $k$ depending on an a {\it priori} information about the
location of $D$ and not the shape. In this paper, we denote by
$\lambda_{j,m}(U)$ for a bounded connected open set $U\subset\Bbb
R^m$ the j-th Dirichlet eigenvalue of $-\triangle$ in $U$.

\subsection{A review of the Colton-Sleeman theorem and Gintides's improvement}

In \cite{CS} Colton and Sleeman have established
the following theorem which we call the Colton-Sleeman theorem.

\proclaim{\noindent Theorem 1.1(\cite{CS}).} Assume that there
exists an open ball $B$ with radius $R$ such that $D\subset B$.  If $k^2<\lambda_{1,m}(B)$,
then $D$ is uniquely determined by $F_D(\,\cdot\,,d;k)$ for
a fixed $d$ and $k$.
\endproclaim

Note that $\lambda_{1,3}(B)=(\pi/R)^2$ and $\lambda_{1,2}(B)=(\gamma_0/R)^2$,
where $\gamma_0$ is the smallest positive zero of the Bessel function
$J_0(z)$.

The assumption means that for fixed $k$, the {\it radius} of $B$
that contains $D$ which is a {\it priori} information about the
location of $D$ can not be large.
It should be pointed out that the optimal case $k^2=\lambda_{1,m}(B)$ is
{\it excluded}.
Their proof does not work for this case.

Their assertion is as follows.

Let $D_1$ and $D_2$ be two obstacles, and $u_1$ and $u_2$ denote the corresponding total fields.
Let $F_1$ and $F_2$ be the corresponding far-field patterns for fixed $d$ and $k$.
If $F_1=F_2$, then $D_1=D_2$.

They employ a contradiction argument which is a traditional one in proving several uniqueness theorems
in inverse obstacle scattering and goes back to Schiffer's idea.

The proof can be divided into five steps.

{\bf\noindent Step 1.}
Assume that the conclusion is not true: $D_1\not=D_2$.

\noindent
{\bf\noindent Step 2.}
Showing $u_1=u_2$ in $D^{\infty}$ with the use of the Rellich lemma \cite{CK},
where $D^{\infty}$ denotes the unbounded connected component of the set
$\Bbb R^m\setminus(\overline D_1\cup \overline D_2)$.

\noindent
{\bf\noindent Step 3.}
Showing the existence of a nonempty connected open set
$D_{\star}\subset D_{\infty}\setminus\overline D_1$ such that $u_1=0$ on
$\partial D_{\star}$, where $D_{\infty}=\Bbb R^m\setminus\overline{D^{\infty}}$
and, if necessary, changing of the role of $u_1$
and $u_2$. See \cite{IS} for this point.

\noindent
{\bf\noindent Step 4.}
Showing $u_1\vert_{D_{\star}}\in H^1_0(D_{\star})$.
This means that there exists a sequence
of smooth functions $\varphi_n\in C^{\infty}_0(D_{\star})$ such that $\varphi_n\longrightarrow u_1$
in $H^1_0(D_{\star})$.
This is because of the general fact: if $U$ is an arbitrary bounded connected open set
and $\varphi\in H^1(U)\cap C^0(\overline U)$ satisfies
$\varphi=0$ on $\partial U$, then $\varphi\in H^1_0(U)$(cf. Corollary 3.28 of \cite{MK}).
Note that the boundary of $D_{\star}$ can be {\it wild} in general and thus one
can not use the characterization of $H^1_0(D_{\star})$ by the trace operator onto $\partial D_{\star}$.
See also \cite{CK, SU} for this point.

\noindent
{\bf\noindent Step 5.}
Showing $k^2\ge \lambda_{1,m}(D_{\star})$.  This is because of $u_1\not\equiv 0$ and $\triangle u_1+k^2 u_1=0$
in $D_{\star}$.

\noindent
{\bf\noindent Step 6.}
Showing $\lambda_{1,m}(D_{\star})\ge\lambda_{1,m}(B)$.
This is because of $D_{\star}\subset B$ and the monotonicity of the first Dirichlet eigenvalue with respect to the domain
which is an implication of the mini-max principle for eigenvalues(the Rayleigh-Ritz formula).

From the last step we have $k^2\ge\lambda_{1,m}(B)$.  Contradiction.

Note that in Step 6 one can not say more like $\lambda_{1,m}(D_{\star})>\lambda_{1,m}(B)$.
This is the reason why the case $k^2=\lambda_{1,m}(B)$ is excluded.

Gintides \cite{G} improved the restriction on $k$ in Theorem 1.1 as
$$\displaystyle
k^2<\lambda_{2,m}(B).
$$

\noindent
His argument after having Step 4 is based on the following four facts.

\noindent
$\bullet$  $\overline{u_1}$ is also satisfies the same Helmholtz equation.

\noindent
$\bullet$  $\overline{u_1}$ and $u_1$ are linearly independent.

\noindent
$\bullet$  the dimension of the first Dirichlet eigenspace is one.  This is the Courant nodal theorem.
See, e.g., on page 133 of \cite{MK} for the proof for a bounded domain {\it without any regularity assumption}
on the boundary just like $D_{\star}$.

\noindent
$\bullet$  the monotonicity of the second Dirichlet eigenvalue with respect to the domain.

\noindent
From these he concludes $k^2\ge \lambda_{2,m}(D_{\star})\ge\lambda_{2,m}(B)$.
Contradiction.

All the argument stated above are based on the {\it multiplicity} of the
eigenvalues and their {\it monotonicity} with
respect to the domain.

{\bf\noindent Remark 1.1.}
Note also that, instead of the monotonicity $\lambda_{1,m}(D_{\star})\ge \lambda_{1,m}(B)$,
Stefanov and Uhlmann in \cite{SU}
used an implication of the Poincar\'e inequality, that is,
$$\displaystyle
\omega_m\le (\lambda_{1,m}(D_{\star}))^{m/2}\text{Vol}\,(D_{\star}),
$$
where $\omega_m$ denotes the volume of the unit ball in $\Bbb R^m$.
They proved a uniqueness theorem at fixed $k$ and $d$ provided $D$ contains a
known obstacle $D_{-}$ and is contained in a known obstacle
$D_{+}$ and $\text{Vol}(D_+\setminus D_{-})<\omega_mk^{-m}$.

\subsection{Statement of the results}

In this paper, we present another method which is based on a real-valued
special function $v$ satisfying $\triangle v+k^2v\le 0$
in a domain that contains unknown obstacles.

Our main result is the following theorem.

\proclaim{\noindent Theorem 1.2.} Let $\Omega$ be an open
connected set with $\overline D\subset\Omega$.
Let $k_0>0$.  Assume that there
exists a real-valued function $v\in C^2(\Omega)$ such that
$\triangle v+k_0^2v\le 0$ in $\Omega$ and $v(x)>0$ for all $x\in\Omega$.
If $k\le k_0$, then $D$ is uniquely determined by
$F_D(\,\cdot\,,d;k)$ for a fixed $d$ and $k$.
\endproclaim

Note that $\Omega$ can be {\it unbounded}; it is assumed that
$\overline D\subset\Omega$ not $D\subset\Omega$.
The $v$ in Theorem 1.2 should be called a {\it supersolution} of
the Helmholtz equation in $\Omega$ at wave number $k_0$ (cf.
\cite{MK} for the notion of the supersolution). Thus Theorem 1.2
can be considered as an application of the supersolution in
inverse obstacle scattering problems. The following corollary
corresponds to Theorem 1.1 including the case when
$k^2=\lambda_{1,m}(\Omega)$.

\proclaim{\noindent Corollary 1.1.}
Let $\Omega$ be a bounded open connected set with $\overline D\subset\Omega$.
If $k^2\le\lambda_{1,m}(\Omega)$,
then $D$ is uniquely determined by
$F_D(\,\cdot\,,d;k)$ for a fixed $d$ and $k$.

\endproclaim
{\it\noindent Proof.} Let $\phi$ be the first Dirichlet eigenfunction
for the negative Laplacian in $\Omega$.  By the Courant nodal theorem, one may assume that
$\phi(x)>0$ for all $x\in\Omega$.  Thus, from Theorem 1.2 with $k_0^2=\lambda_{1,m}(\Omega)$
and $v=\phi$, one obtains the desired uniqueness result.

\noindent
$\Box$

However, this fact itself is {\it not new} since the result
is a special case of the result by Gintides \cite{G} as mentioned in Subsection 1.1
under the condition $k^2<\lambda_{2,m}(B)$ and $D\subset \Omega\equiv B$.

{\bf\noindent Example 1.}
Let $B=\{x\in\Bbb R^m\,\vert\,\vert x\vert<R\}$.
For $\Omega=B$ one can choose
$$\displaystyle
\phi(x)=\left\{
\begin{array}{lr}
\displaystyle
J_0(k_0\vert x\vert), & \quad\text{if
$m=2$,}\\
\\
\displaystyle
\frac{\sin\,k_0\vert x\vert}{\vert x\vert}, & \quad\text{if $m=3$,}
\end{array}
\right.
$$
where $k_0=\lambda_{1,m}(B)$.

When $\Omega$ is bounded, one can not find a positive supersolution of the Helmholtz equation
in $\Omega$ at the wave number $k>\sqrt{\lambda_{1,m}(\Omega)}$.  This is because
of the following fact.

\proclaim{\noindent Proposition 1.1.} Let $\Omega$ be a bounded open
connected set of $\Bbb R^m$. There exists a real-valued function
$v\in C^2(\Omega)$ such that $\displaystyle\triangle v+k^2v\le 0$
in $\Omega$ and $v(x)>0$ for all $x\in\Omega$ if and only if
$k^2\le\lambda_{1,m}(\Omega)$.
\endproclaim

For the proof see Appendix. Thus Theorem 1.2 does not yield
a new result beyond the Colton-Sleeman theorem and Gintides's result
when $\Omega$ is bounded.
However, when $\Omega$ is {\it
unbounded}, there is a possibility of having a positive
supersolution in $\Omega$. This is an advantage of Theorem 1.2.
The following two corollaries are new and not covered by the
Colton-Sleeman theorem or Gintides's result.

\proclaim{\noindent Corollary 1.2.}
Let $\Omega'$ be a bounded open connected set of $\Bbb R^2$ with
$\overline D\subset\Bbb R\times\Omega'$.
If $k^2\le\lambda_{1,2}(\Omega')$,
then $D$ is uniquely determined by
$F_D(\,\cdot\,,d;k)$ for a fixed $d$ and $k$.
\endproclaim

\noindent
{\it Proof.} Let $\phi'$ be the first positive Dirichlet eigenfunction
for the negative Laplacian in $\Omega'$.
Define $v(x_1,x_2,x_3)=\phi'(x_2,x_3)$ for $x\in\Bbb R\times\,\Omega'$.
This $v$ satisfies $\triangle v+k_0^2v=0$ in $\Bbb R\times\,\Omega'$ with $k_0^2=\lambda_{1,2}(\Omega')$
and $v(x)>0$ for all $x\in\Bbb R\times\,\Omega'$.

\noindent
$\Box$

{\bf\noindent Example 2.}  Let $\Omega'=]-R\,,R[\times\,]-h\, ,h[$ with $h, R>0$.
Then
$$\displaystyle
\lambda_{1,2}(\Omega')
=\left(\frac{\pi}{2}\right)^2\left(\frac{1}{h^2}+\frac{1}{R^2}\right)
$$
and
$$\displaystyle
\phi'(x_2,x_3)=\cos\,\frac{\pi}{2R}x_2\,
\cos\,\frac{\pi}{2h}\,x_3,\,\,(x_2,x_3)\in\,]-R\,,R[\times\,]-h\, ,h[.
$$
Thus the condition $k^2\le\lambda_{1,2}(\Omega')$ is equivalent to
$$\displaystyle
k\le\frac{\pi}{2}\sqrt{\frac{1}{h^2}+\frac{1}{R^2}}.
$$

A similar idea yields

\proclaim{\noindent Corollary 1.3.}
Let $J$ be a bounded open interval of $\Bbb R$ with
$\overline D\subset\Bbb R^2\times J$.
If $k^2\le\lambda_{1,1}(J)$,
then $D$ is uniquely determined by
$F_D(\,\cdot\,,d;k)$ for a fixed $d$ and $k$.
\endproclaim

{\bf\noindent Example 3.}
Let $J=]-h,\,h[$ with $h>0$.
Then $\lambda_{1,1}(J)=(\pi/2h)^2$ and an associated positive Dirichlet eigenfunction for the negative Laplacian in $J$ is given by
$$\displaystyle
\phi(x_3)=\cos\,\frac{\pi}{2h}x_3,\,\,\vert x_3\vert<h.
$$
The condition $k^2\le \lambda_{1,1}(J)$ is equivalent to
$$\displaystyle
k\le\frac{\pi}{2h}.
$$

Examples 2 and 3 suggest that the larger is the number of unbounded directions of the domain $\Omega$,
the lower is the upper bound $k_0$.

\section{Proof of Theorem 1.2}

We start with describing a well known identity.

\proclaim{\noindent Proposition 2.1.} Let $u$ and $v$ be arbitrary
smooth functions on an open set $U$ and satisfy $v(x)\not=0$ for
all $x\in U$.  Then we have
$$\displaystyle
\nabla\cdot(v^2\nabla\varphi)
=v\triangle u-u\triangle v\,\,\text{in}\,U,
$$
where
$$\displaystyle
\varphi=\frac{u}{v}.
$$

\endproclaim

Using this identity, we have the following lemma.

\proclaim{\noindent Lemma 2.1.}
Let $k_0>0$.  Assume that there
exists a real-valued function $v\in C^2(\Omega)$ such that
$\triangle v+k_0^2v\le 0$ in $\Omega$ and $v(x)>0$ for all $x\in\Omega$.
Let $U$ be a bounded open connected set of $\Omega$ with $\overline U\subset\Omega$.
If $k\le k_0$ and $u\in C^2(U)\cap C^0(\overline U)$ satisfies
$$\displaystyle
\begin{array}{c}
\displaystyle
\triangle u+k^2u=0\,\,\text{in}\,U,\\
\\
\displaystyle
u=0\,\,\text{on}\,\partial U,
\end{array}
$$
then $u=0$ in $U$.
\endproclaim
{\it\noindent Proof.}
Define $\varphi=u/v$ in $U$.
It follows from Proposition 2.1 that
$$\begin{array}{c}
\displaystyle
\nabla\cdot(v^2\nabla\varphi)+(\triangle v+k^2v)v\varphi=0\,\,\text{in}\, U.
\end{array}
\tag {2.1}
$$
Since $v^2$ has a {\it positive uniform lower bound} on $U$, $\varphi=0$ on $\partial U$ and
$$\displaystyle
(\triangle v+k^2v)v=(\triangle v+k_0^2v)v-(k_0^2-k^2)v^2\le 0\,\,\text{in}\,U,
$$
the {\it weak maximum principle} (\cite{MK}) yields
$\varphi=0$ in $U$ and thus $u=0$ in $U$.

\noindent
$\Box$

The proof of Theorem 1.2 starts with having Step 3.
Applying Lemma 2.1 with $u=u_1$ and $U=D_{\star}$, we have $u_1=0$ in $D_{\star}$.
Then the unique continuation gives $u_1=0$ in $\Bbb R^m\setminus\overline D_1$ and this contradicts
$u_1\sim e^{ikx\cdot d}$ as $\vert x\vert\longrightarrow\infty$.
Therefore it must hold that $D_1=D_2$.

\noindent
$\Box$

{\bf\noindent Remark 2.1.}
Note that: if one starts with having Step 4 in Subsection 1.1, then $\varphi=u_1/v\vert_{D_{\star}}\in H^1_0(D_{\star})$.
Using (2.1) and a sequence in $C_0^{\infty}(D_{\star})$ that converges to $\varphi$
in $H^1_0(D_{\star})$, we have
$$\displaystyle
\int_{D_{\star}}v^2\vert\nabla\varphi\vert^2dx-\int_{D_{\star}}(\triangle v+k^2)v\vert\varphi\vert^2dx=0.
$$
This also yields the same conclusion as above.
This avoids the use of the weak maximum principle, however, needs a knowledge that $u_1\vert_{D_{\star}}\in H^1_0(D_{\star})$.

{\bf\noindent Remark 2.2.}
The argument done in the proof of Lemma 2.1 together with the use of  Proposition 2.1
are well-known typical one in studying the maximum principle for general elliptic partial differential equations
(e.g., \cite{MK} and introduction of \cite{BNV}).  Here we presented it just for the use of (2.1) in Remark 2.1,
i.e., its {\it use} in inverse obstacle scattering problems.

\section{Conclusion}

The previous known proof of the Colton-Sleeman and Gintides's
improvement are based on some facts on the Dirichlet eigenvalues
of the negative Laplacian in a domain and their monotone
dependence on domains.

Our proof is extremely simple and uses a positive {\it supersolution} $v$ of
the Helmholtz equation in a domain $\Omega$ that contains the closure of
unknown obstacle. Domain $\Omega$ in three dimensions can be
{\it unbounded} for a {\it single direction} as shown in Corollary 1.2
and {\it two directions} as in Corollary 1.3 
if the wave number has a bound depending on the size of the ``bounded part''
of $\Omega$.

$$\quad$$

\centerline{{\bf Acknowledgement}}

This research was partially supported by Grant-in-Aid for
Scientific Research (C)(No. 21540162) of Japan  Society for the
Promotion of Science.
The author would like to thank a referee who made constructive suggestions including 
the change of the title of the paper.

\section{Appendix. Proof of Proposition 1.1}

First assume the existence of $v$ and that the
conclusion is not true. Thus one has $k^2>\lambda_{1,m}(\Omega)$.  Let
$\Omega_1, \Omega_2,\cdots$ be an exhaustion of $\Omega$ from
below in the sense that each $\Omega_j$ is open connected and
$\overline\Omega_{j-1}\subset\Omega_{j}\uparrow\Omega$. Then we
have $\lambda_{1,m}(\Omega_j)\ge\lambda_{1,m}(\Omega)$ and
$\lambda_{1,m}(\Omega_j)\longrightarrow\lambda_{1,m}(\Omega)$.
The latter is also a well known consequence of the Rayleigh-Ritz formula.
Thus for a large $j$ we have $k^2>\lambda_{1,m}(\Omega_j)$. Let
$\phi_j$ denote the first Dirichlet eigenfunction for the negative
Laplacian in $\Omega_j$. Define $\varphi=\phi_j/v$ in $\Omega_j$.
Since $\overline{\Omega_j}\subset\Omega$, $\varphi$ belongs to
$C^2(\Omega_j)\cap C^0(\overline\Omega_j)$. It follows from
Proposition 2.1 that
$$\begin{array}{c}
\displaystyle
\nabla\cdot(v^2\nabla\varphi)+(\triangle v+k_0^2v)v\varphi=0\,\,\text{in}\,\Omega_j,
\end{array}
$$
where $k_0^2=\lambda_{1,m}(\Omega_j)$.  Since $\varphi=0$ on
$\partial\Omega_j$ and $(\triangle v+k_0^2v)v=(\triangle
v+k^2v)v+(k_0^2-k^2)v^2\le 0$ in $\Omega_j$, the maximum principle
yields $\varphi=0$ in $\Omega_j$; however, this is impossible by
the Courant nodal theorem. Therefore it must hold that
$k^2\le\lambda_{1,m}(\Omega)$.  Conversely if
$k^2\le\lambda_{1,m}(\Omega)$, then choose the first positive Dirichlet
eigenfunction $\phi$ for the negative Laplacian in $\Omega$ and
set $v=\phi$.  Then $\displaystyle\triangle
v+k^2v=\triangle\phi+k_1^2\phi +(k^2-k_1^2)\phi\le 0$ in $\Omega$,
where $k_1^2=\lambda_{1,m}(\Omega)$.

\noindent
$\Box$

\vskip1cm
\noindent
e-mail address

ikehata@math.sci.gunma-u.ac.jp


\begin{thebibliography}{99}



\bibitem{BNV}  Berestycki, H., Nirenberg, L. and Varadhan, S. R. S.,
               The principle eigenvalue and maximum principle for second-order
               elliptic operators in general domains,
               Comm. Pure Appl. Math., {\bf 47}(1994), 47-92.


\bibitem{CS}  Colton, D. and Sleeman, B. D.,
              Uniqueness theorems for the inverse problem of acoustic scattering,
              IMA J. Appl. Math., {\bf 31}(1983), 253-259.


\bibitem{CK} Colton, D. and Kress, R., Inverse acoustic and electromagnetic scattering theory,
             Springer, second edition, 1998.



\bibitem{G}  Gintides, D., Local uniqueness for the inverse scattering problem in acoustics
             via the Faber-Krahn inequality, Inverse Problems, {\bf 21}(2005), 1195-1205.




\bibitem{IS} Isakov, V., Inverse problems for partial differential equations,
             Springer, second edition, New York, 2006.



\bibitem{MK} Murata, M. and Kurata, K., Elliptic and parabolic type partial differential equations,
             in japanese, Iwanami Shoten, 2006, Tokyo.




\bibitem{SU} Stefanov, P. and Uhlmann, G.,
             Local uniqueness for the fixed energy fixed angle inverse problem in obstacle scattering,
             Proc. AMS., {\bf 132}(2004), 1352-1354.







\end{thebibliography}
\end{document}